\documentclass[12pt, a4paper]{article}
\usepackage[utf8]{inputenc}
\usepackage{amscd,amssymb,amsthm,amsmath,comment}
\title{A remark on cohomology of nilpotent Lie algebras}
\author{Grigory Papayanov}
\date{2023}

\newcommand{\g}{{\mathfrak g}}

\renewcommand{\k}{{\Bbbk}}

\newcommand{\arrow}{{\:\longrightarrow\:}}
\newcommand{\Z}{{\mathbb Z}}

\newcommand{\R}{{\mathbb R}}

\def\1{\sqrt{-1}\:}

\renewcommand{\phi}{\varphi}
\renewcommand{\epsilon}{\varepsilon}

\newcommand{\Sym}{\operatorname{Sym}}

\newcommand{\Id}{\operatorname{Id}}

\newcommand{\cinf}{C_{\infty}}

\newcommand{\gr}{\operatorname{gr}}

\renewcommand{\Im}{\operatorname{Im}}

\newcommand{\black}{\bullet}

\newcommand{\ainf}{A_{\infty}}

\newcommand{\Cobar}{\Omega}
\renewcommand{\Bar}{\mathrm{B}}
\newcommand{\colim}{\operatorname{colim}}
\newcommand{\Ext}{\operatorname{Ext}}
\newcommand{\FREE}{FREE}
\newcommand{\Har}{\mathrm{T}}
\newcommand{\Tor}{\operatorname{Tor}}

\newcounter{Mycounter}[section]
\newcounter{lemma}[section]
\renewcommand{\thelemma}{\noindent{Lemma \thesection.\arabic{lemma}}}
\newcommand{\lemma}{%
     \setcounter{lemma}{\value{Mycounter}}
     \refstepcounter{lemma}
     \stepcounter{Mycounter}
     {\bf \thelemma:\ }}

\newcounter{claim}[section]
\newcounter{sublemma}[section]
\newcounter{corollary}[section]
\newcounter{theorem}[section]
\renewcommand{\thetheorem}{\noindent{Theorem \thesection.\arabic{theorem}}}
\newcommand{\theorem}{%
     \setcounter{theorem}{\value{Mycounter}}
     \refstepcounter{theorem}
     \stepcounter{Mycounter}
     {\bf \thetheorem:\ }}	

\newcounter{conjecture}[section]
\newcounter{proposition}[section]
\newcounter{definition}[section]
\renewcommand{\thedefinition}
       {\noindent{Definition~\thesection.\arabic{definition}}}
\newcommand{\definition}{%
     \setcounter{definition}{\value{Mycounter}}
     \refstepcounter{definition}
     \stepcounter{Mycounter}
     {\bf \thedefinition:\ }}

\newcounter{example}[section]
\newcounter{remark}[section]
\renewcommand{\theremark}{\noindent{Remark \thesection.\arabic{remark}}}
\newcommand{\remark}{%
     \setcounter{remark}{\value{Mycounter}}
     \refstepcounter{remark}
     \stepcounter{Mycounter}
     {\bf \theremark:\ }}

\newcounter{problem}[section]
\newcounter{question}[section]
\begin{document}

\begin{center}
{\Large\bf
A remark on cohomology of nilpotent Lie algebras
}
\\[4mm]

Grigory Papayanov
\\[6mm]

\end{center}

\begin{center}
\textbf{Abstract}
We prove that the cohomology algebra of a
conilpotent Lie coalgebra is generated in
degree 1 as an $\ainf$-algebra. By dualizing,
the same is true about cohomology of finite dimensional nilpotent
Lie algebras. In the process, we provide a proof of the conilpotent
version of the dual Poincare-Birkhoff-Witt theorem.
\end{center}

\section{Introduction}

This note grew out of desire to calculate the cohomology
of free nilpotent Lie algebras $L^{\le n}(V) := L(V)/L^{>n}(V)$.
Here $V$ is a vector space, $L(V)$ is the free Lie algebra
generated by $V$, naturally graded by the commutator length:
$L(V) = \oplus_{n \ge 1} L^n(V)$, and $L^{\le n}(V)$ and $L^{>n}(V)$
are, respectively, the quotient algebras and the 
ideals associated to this grading.
Since a free nilpotent Lie algebra is a very natural object, it
seems that its cohomology should also be a very natural object,
and hence any understanding of its structure might be desirable.
Indeed, even in the simplest nontrivial case $n = 2$ the calculation
of the $GL(\R^m)$-character of the representation $H(L^{\le 2}(\R^m))$
performed by J\'ozefiak-Weyman (\cite{JOZWEY}) and by
Sigg (\cite{SIGG}) gives a categorification of the well-known
Littlewood identity:

$$
\prod_{i=1}^m (1-x_i) \prod_{1 \le i < j \le m} (1-x_ix_j) = 
\sum_{\lambda} (-1)^{\frac{1}{2}(|\lambda|+r(\lambda))}s_\lambda(x).
$$

Here the sum in the right-hand part is the sum of the Schur polynomials
corresponding to all self-conjugate Young diagrams. By $|\lambda|$
is denoted the number of all squares in the diagram and by $r(\lambda)$ the
number of diagonal squares.

\hfill

So it is plausible that knowledge of cohomology of free nilpotent Lie 
algebras of higher nilpotency degree may bring some insights
on other mathematical questions --- at the very least, it will be
a source of new identities on symmetric functions. However, for now
 only partial results are available,
 mostly obtained by direct
calculations, see, e.g., \cite{TIRAO}. 

\hfill

It is often the case that the existence of extra structures 
helps calculations by limiting their scope. In the case of
cohomology calculations one such extra structure is the multiplication.
Moreover, the multiplication on cohomology 
is part of a richer structure --- the $\cinf$ structure (\cite{KADEISHC}),
which consists
of a collection of of polylinear operations
that refine Massey products in a certain way.
For semisimple $\g$, the algebra $H(\g)$
knows little information about $\g$ itself ---
it is just the exterior algebra on some generators
 readable from the Dynkin diagram.
There is no nontrivial $\cinf$-structure in this case.
Nilpotent algebras are the opposite: knowing the $\cinf$-structure, one 
can reconstruct from it the original Lie algebra up to isomorphism.
Using the calculations of Sigg, Dubois-Violette and Popov in \cite{DVP}
were able to prove that as a $\cinf$-algebra, the cohomology algebra
$H(L^{\le 2}(V))$ is generated by $H^1$ with only $m_2$ and $m_3$ operations.

\hfill

Since free nilpotent Lie algebras admit such a simple
(four words!) description, and since their cohomology algebras
are tied very closely to them, one is lead to think that perhaps
 $\cinf$-algebra 
$H(L^{\le n}(V))$ also admits a short description. From the results
of this paper it follows that as a $\cinf$-algebra it is generated by its 
component of degree 1. However, not much is clear about the relations.
It is not hard to calculate $H^2$ --- indeed, the exact sequence

$$
0 \arrow L^{n+1}(V) \arrow L^{\le n+1}(V) \arrow L^{\le n}(V) \arrow 0
$$
 
\hfill

is the universal central extension of $L^{\le n}(V)$, hence 
$H^2(L^{\le n}(V)) = (L^{n+1}(V))^\vee$. The $\cinf$-operations
$m_i: H^1 = L^i(V^\vee) \arrow (L^{n+1}(V))^\vee$ are either zero
when $i \ne n+1$ or natural isomorphisms when $i = n+1$. However,
the behaviour in higher degrees remains elusive. To have
an example of complications arising when trying to generate infinity-algebras
by generators and relations, the reader may consult \cite{FOUR}.

\hfill

On the positive side, the result about 1-generatedness
is not particular to the case of free nilpotent Lie algebras:
the results of this paper show that the cohomology algebra $H(\g)$
is $\cinf$-generated by $H^1(\g)$ for {\em any nilpotent} finite dimensional
Lie algebra $\g$. It also follows that for inverse limits
of finite dimensional nilpotent Lie algebras the $\cinf$-algebra
of {\em continious cohomology} $H_c(\g)$ is generated by $H^1$. This result
is probably known to some experts in the field
 (especially in the case where there exists a positive grading on $\g$ --- 
the general case of nilpotent $\g$ seems to be a more obscure knowledge).
However, since from time to time papers with partial results appear
(see, for example, \cite{MILLION} dealing with the case of the Witt algebra),
there might be a demand for a paper where 1-generatedness result
is written.

\hfill

We prove the result for conilpotent Lie coalgebras rather than nilpotent
Lie algebras, because the proof is more natural 
and we feel that it is more easily digestible that way. First, this step
from algebras to coalgebras removes the necessity to have a positive
grading. Second, in this way we don't have to mention the
cobar construction at all, and it seems
that mathematicians are in general more acquainted
with the bar construction. Even though our proof is just an 
application of the general cobar-bar duality (and we formulate
a theorem for a general pair of Koszul operads), we wanted to make
the paper largely self-contained --- 
the inclusion of the discussion of Koszul duality
for operads would probably double its length.
 This is also the reason we're working
with $\ainf$-algebras and not $\cinf$-algebras: the former
type algebras are familiar to more mathematicians, and the latter
is sometimes perceived as exotic. As a side consequence of our choice
of approach, we obtain a proof of the analogue of the Poincare-Birkhoff-Witt
theorem for the conilpotent coenveloping coalgebra of a conilpotent Lie algebra.
We claim no originality for our proof, since it is just the dualization
of the Positselski's proof of the PBW theorem (\cite{POSIC}). However, since
Positselski wrote down his proof in great generality --- for 
nonhomogeneous quadratic deformations of Koszul algebras --- it might be
useful to someone to have it in the traditional formulation.

\hfill

Our proof works only over the field of characteristic zero.
The only time the characteristic comes up is during the 
calculation of $\Ext_{\Lambda(V)}(\k, \k)$, but
this calculation is crucial to the proof.
Reader, interested in the cohomology of
free nilpotent Lie algebras over $\Z$, may find
\cite{ROMANO} useful.

\hfill

{\bf Acknowledgments:}
The author is grateful to Kostya Tolmachov, Sasha Viktorova
and Sasha Zakharov
for their interest in the work and for their helpful comments.
The author especially wishes to thank Volodya Dotsenko for telling him
that this result is worth writing down. This paper was written
during a visit to the Weizmann Institute of Sciences. I am very
grateful for its hospitality.

\hfill

\section{Lie coalgebras}

\definition
A {\bf Lie coalgebra} over a field $\k$ is a vector
space $L$ together with the comultiplication operation
$\Delta: L \arrow L \otimes L$ satisfying the skew-cocommutativity and
coJacobi relations. The skew-cocommutativity relation
means that the image of $\Delta$ lies in the subspace 
$\Lambda^2L \subset L \otimes L$ of skew-symmetric tensors. The coJacobi
relations could be phrased as follows: let $\xi$ be an
operation of cyclic permutation on $L \otimes L \otimes L$.
Then the coJacobi relation tells that the following map 
$L \arrow L \otimes L \otimes L$

\begin{equation}
(\Id + \xi + \xi^2) \circ (\Id \otimes \Delta) \circ \Delta
\end{equation}

must vanish.

\hfill

\remark
The dual space $L^\vee$ to a Lie coalgebra $L$ is a Lie algebra,
the converse, verbatim, is not true, but it becomes true if one 
considers the space $L^\vee$ as a {\em topological} vector
space with its natural Tychonoff topology. Topological vector spaces
which are dual to discrete vector spaces are precisely
{\em linearly compact} spaces (see, e.g., \cite{DIEU}),
and the space $W^\vee$ of continious linear functionals on a linearly
compact vector space $W$ is precisely the vector space of which
$W$ is the dual. As a result, the categories of Lie coalgebras
and of linearly compact Lie algebras are (anti-)equivalent.

\hfill

\definition
A {\bf Chevalley-Eilenberg} algebra $C(L)$ of a Lie coalgebra $L$ is
the exterior algebra $\Lambda L$ with the grading
$C^k(L) := \Lambda^kL$ and 
the differential $\delta$ which on the generator
space $L$ is equal to $\Delta$. It is naturally augmented
by the projecton on the degree zero component.

\hfill

\remark
To obtain $C(L)$ from $L^\vee$ one has to take the 
space of continious polylinear functionals on $L^\vee$.

\hfill

The calculation which can be found in, e.g., \cite{MICHAELIS},
shows that $\delta^2$ is zero and hence $C(L)$ is a dg-algebra.

\hfill

\definition\label{conilpotent}
A coalgebra $L$ is called {\bf filtered} or {\bf conilpotent}
if its dual topological Lie algebra is a projective limit of finite
dimensional nilpotent
Lie algebras. Alternativelty, it means
that there exists an increasing ascending exhaustive filtration
by subspaces $F_*L$
such that $F_0L = 0$ and such that $\Delta(F_i)$ lies
inside $\sum_{p+q=i} F_{p} \otimes F_{q} \subset L \otimes L$.

\hfill

\remark
The previous definition could be applied to other kinds of coalgebras,
not necessarily Lie ones. Beware, though, that in the coaugmented
case one may require $F_0L$ to be equal to the image of coaugmentation.
In this article all coalgebras are not supposed to have neither a counit
nor a coaugmentation.

\hfill

The tensor powers of conilpotent Lie coalgebra $L$
inherit the filtration and hence its Chevalley-Eilenberg
algebra $C = C(L)$ is then also endowed with an ascending
exhaustive filtration such that $F_0C = \k = C^0$. The
Chevalley-Eilenberg differential preserves this filtration.

\hfill

For a linearly compact Lie algebra $\g$ we will denote by $U\g$ its 
topological universal enveloping associative algebra. 
It is a quotient of the completed (both in the topology
of $\g$ and in tensor degree topology) free algebra generated by $\g$
by the closure of the ideal generated by elements
of the form $xy - yx = [x, y]$. This algebra is naturally
augmented. We will denote its augmentation ideal by $I$.

\hfill

\definition
Let $L$ be a Lie coalgebra, and let $\g := L^\vee$
be the corresponding dual linearly compact Lie algebra. 
Consider the coassociative coalgebra

\begin{equation}
UL := \colim_k (U\g/I^k)^\vee.
\end{equation}

This coalgebra is called a {\bf conilpotent coenveloping coalgebra} of $L$.

\hfill

In contrast to the
algebra case, there exist Lie coalgebras for which
conilpotent coenveloping coalgebras
are trivial. One can take, for example,
the Lie coalgebra dual to any semisimple Lie algebra.
 However, for conilpotent Lie
coalgebras an exact analogue of the Poincare-Birkhoff-Witt theorem holds.
One proof of this theorem is available in the appendix
of \cite{BETTS-DOGRA}. We will provide another one, which
basically coincides with the proof of the usual PBW theorem
due to Positselski (\cite{POSIC}).

\hfill

\section{A-infinity algebras}

\definition
Let $A$ be an associative dg-algebra, not necessarily unital.
 Its {\bf bar-construction}
is a coassociative dg-coalgebra $\Bar A$ which, as a coalgebra,
is the tensor coalgebra generated by $A[1]$, and whose codifferential
is the unique codderivation of the tensor coalgebra
whose corestriction to $A[1]$ is given by the differential and multiplication
in $A$. Modulo signs, we have

\hfill

\begin{equation}\label{bard}
d_{\Bar}(a_1 \otimes \dots \otimes a_n) = 
\sum \pm a_1 \otimes \dots da_k \dots \otimes a_n + 
\sum \pm a_1 \otimes \dots a_ia_{i+1} \dots \otimes a_n.
\end{equation}
 
\hfill

The condition that $d_{\Bar} = 0$ is precisely the defining
equations for an associative dg-algebra: the multiplication
is associative, the differential is square-zero and satisfies
the Leibniz rule. The bar-construction is naturally 
(moreover, uniquely) coaugmented
by the inclusion of tensors of degree zero.
 A morphism $f: A \arrow B$ induces
a map of dg-coalgebras, which we will denote by the
same letter $f: \Bar A \arrow \Bar B$. This map preserves the coaugmentations
of bar-constructions.

\hfill

\definition
An {\bf A-infinity} algebra $A$ is a graded vector space
endowed with the square-zero coderivation $d_\Bar A$ 
of the tensor coalgebra generated by $A[1]$ preserving 
its unique coaugmentation. We will denote the
resulting dg-coalgebra by $\Bar A$ also. A {\bf morphism}
between two $\ainf$-algebras $A$ and $B$ is, by definition,
a morphism of dg-coalgebras $\Bar A \arrow \Bar B$. It
preserves the coaugmentations automatically.

\hfill

From the universal property of the tensor coalgebra
it follows that the codifferential $d_\Bar$ is determined
by its corestriction to $A[1]$, that is, by a collection of maps

\begin{equation}
m_i: A[1]^{\otimes i} \arrow A[2],
\end{equation}

satisfying a collection of quadratic equations corresponding
to $d_\Bar^2=0$
and an $\ainf$-morphism $f: \Bar A \arrow \Bar B$
is determined by a collection of maps

\begin{equation}
f_i: A[1]^{\otimes i} \arrow B[1],
\end{equation}

satisfying a collection of quadratic equations corresponding
to $d_{\Bar B}f = fd_{\Bar A}.$ It is clear that dg-algebras are
precisely those $\ainf$-algebras for which $m_i = 0$ for $i \ne 1,2$.
However, a notion of an $\ainf$-morphism is more general than the notion
of a morphism of dg-algebras, except for the case
when both algebras are concentrated in a single degree $\ne 1$. We will
call the collections $\{m_i\}$ and $\{f_i\}$ the
{\bf Taylor components} of the entities defined by them.

\hfill

Note that for a dg-algebra $A$ the bar-construction $\Bar A$ is a
filtered coalgebra. The filtration is given by 

\begin{equation}
F_i\Bar A := \bigoplus_{k=0}^i (A[1])^{\otimes k}.
\end{equation}

If an $\ainf$-algebra $A$ is determined by Taylor components
$(m_1, m_2, \dots)$, then $\gr_F\Bar A$ is $\oplus_{i} A[1]^{\otimes i}$
with the differential $m_1$. In particular, $m_1^2 = 0$, and it defines
a structure of a complex on $A$. 

\hfill

Note that if $f: \Bar A \arrow \Bar B$ is an $\ainf$-morphism,
then it automatically preserves this filtration. If the Taylor components
of $f$ are $(f_1, f_2, \dots)$, then $\gr_Ff = f_1$. In particular,
$f_1m_{A,1} = m_{B,1}f_1$.

\hfill

\definition
An $\ainf$-morphism $f: \Bar A \arrow \Bar B$ is called a quasiisomorphism
if it is a filtered quasiisomorphism of dg-coalgebras, that is,
if $\gr_Ff$ is a quasi-isomorphism. Equivalently, since 
tensor products are exact over a field, $f_1$ is a quasiisomorphism.

\hfill

\remark\label{ssequence}
It follows then that $f$ is a quasi-isomorphism of dg-coalgebras,
if $f$ is an $\ainf$-quasiisomorphism. Indeed, 
from 5-lemma it follows that
$f \vert_{F_k}$ is a quasi-isomorphism for any $k$. Since
$\Bar = \colim_k F_k \Bar$ the result follows from the exactness
of the direct limit.

\hfill

The usefulness of $\ainf$-algebras is illustrated by the following
theorem, which is usually attributed to Kadeishvili:

\hfill

\definition
An $\ainf$-algebra $A$ is called {\bf minimal}
if $\gr_F\Bar A$ has zero differential.

\hfill

\theorem
Let $(A, m_1, m_2, \dots)$ be an $\ainf$ algebra. Let $H$ be the cohomology 
of the differential $m_1$. Then $H$  admits a structure of a minimal
$\ainf$-algebra and there is an $\ainf$ quasi-isomorphism
between this structure and $A$. Moreover, the structure on $H$ is unique
up to an $\ainf$-{\em isomophism}.

\hfill

\remark
The theorems of this kind were considered by rational homotopists
before Kadeishvili: see, for example, \cite[Theorem 2.2]{SULLIVAN}.
However, the actual Kadeishvili's result (\cite{KADEISHA})
is more general than the theorem above: for example,
he doesn't require the base ring to be a field,
but only requires $H$ to be a projective module. See \cite{PETERSEN}
for discussion.

\hfill

\section{1-generation}

Let $(A, m_2, m_3, \dots)$ be a minimal $\ainf$-algebra.
We are interested in the cohomology of its bar-complex $H(\Bar A)$.
As $\Bar A$ is a filtered (by the tensor filtration) complex, 
its cohomology $H(\Bar A)$ inherit the filtration.
That is, $F^kH(\Bar A)$ are classes that could be represented by
cocycles that lie in $F^k\Bar A$. 

\hfill

\definition
A minimal $\ainf$-algebra $A$ is called 
{\bf 1-generated} if it is positively graded and
$F^1H^j(\Bar A) = 0$ for $j \ge 1$.

\hfill

\remark
Note that this notion is well-defined in the sense
that it is invariant under $\ainf$-automorphisms.
Indeed, an $\ainf$-automorphism of $A$
defines a filtered automorphism
of the dg-coalgebra $\Bar A$.

\hfill

The following lemma is a generalization to the
$\ainf$ case of a well-known criterion of 1-generatedness
of graded algebras, see e.g. \cite{POSVIS}

\hfill

\lemma
A positively graded $\ainf$ algebra $A$ is 1-generated if and only if
any element $x \in A^k$ can be expressed as 
a linear combination of iterated
compositions of maps of the form 
$\Id^{\otimes i} \otimes m_j \otimes \Id^{\otimes k}$
applied to elements in $A^{\otimes n}$.

\hfill

\remark
This is equivalent to the fact that the map
$\FREE(A^1) \arrow A$, where $\FREE$ is the functor
of the free minimal $\ainf$-algebra
(which exists since minimal $\ainf$-algebras
are algebras over a certain operad, see e.g. \cite{LODAYVALETTE}), is surjective.
The fact that this is invariant under $\ainf$-isomorphisms
is less obvious then the previous definition of 1-generatedness.

\hfill

{\bf Proof:}
First remember that the space $F^1H^k(\Bar A)$ is the subspace
of $H^k(\Bar A)$ that could be represented by classes 
in $F^1(\Bar A) = A[1] \subset \Bar (A)$. Since $A$ is minimal,
the subspace $A[1]$ consists of closed elements, hence
$F^1H^k(\Bar A)$ is the quotient of $A^{k-1}$ by the image of
 $d_{\Bar A}$.
Now we can prove the equivalence of two definitions by induction.
Note that $\Im d_{\Bar A} \cap A^2$ is spanned by images of arrows
$(A^1)^{\otimes i} \arrow A^2$ which up to sign coincide with the
operations $m_i$. Hence $A^2$ is generated from $m_i$'s by $A^1$
if and only if $F^1H^1(\Bar A) = 0$. Next, suppose we proved the equivalence
for all $A^k$, where $k < n$ for some $n$. Then space  
$\Im d_{\Bar A} \cap A^n$ is spanned by images of $m_i$'s applied to
the subspaces $\otimes_{j=2}^i A^{p_j}$, where 
$\sum p_j = n-i+2$. Again, by degree reasons all $p_j$'s are less then $n$,
and the vanishing of $A^n/\Im d_{\Bar A}$ is equivalent to the fact
that the elements of $A^n$ are expressable via
the operations $m_i$ through the elements of smaller degree, hence, by the 
base of induction, they are expressable through elements of $A^1$.

\hfill

\section{Conilpotent PBW theorem}\label{PBW}

The PBW theorem will follow from the following computation
of the bar-construction of the Chevalley-Eilenberg algebra
of a conilpotent Lie coalgebra $L$. Denote by $C^{>0}(L)$
the quotient of $C(L)$ by the coaugmentation map $\k \arrow C(L)$.

\hfill

\lemma\label{barcobar}
For a conilpotent Lie coalgebra $L$, the
cohomology of $\Bar C^{>0}(L)$ vanish in positive degrees.
The coalgebra $H^0(\Bar C^{>0}(L))$ is isomorphic to $U(L)$,
the conilpotent coenveloping coalgebra of $L$.

\hfill

{\bf Proof:} We have two naturally defined filtrations
on $\Bar C^{>0}(L)$. The "stupid" filtration on $C^{>0}(L)$
given by $G^k = \oplus_{i \ge k} C^{i}(L)$
extends to a filtration on $\Bar C^{>0}(L)$, which we will denote by
  $G$. The filtration induced
from the conilpotent filtration on $L$ will be denoted by $N$. 
The filtration $G$ is descending and non-complete,
the filtration $N$ is ascending and exhaustive.

\hfill

Consider the complex $\gr_N\Bar C^{>0}(L)$. The filtration $G$
induces a filtration on it, which we will also denote by $G$. An important
fact is that $G$ is finite on each $N_{i}/N_{i-1}$. Consider now the
complex $\gr_G\gr_N \Bar C_{>0}(L)$. We want to show that its higher cohomology
vanish. From finiteness of $F$ and exhaustiveness of $N$
it will follow that higher cohomology of $\Bar C^{>0}(L)$ would
vanish as well.

\hfill

Note that the double associated graded
$\gr_G\gr_N \Bar C^{>0}(L)$ does not depend on the cobracket
on $L$ and only depends on the underlying vector space.
Indeed, $\Bar C^{>0}(L)$ is naturally a bicomplex, namely,
$\Bar C^{>0}(L) = \oplus_{i, j \ge 0} B^{i,j}$, where

\begin{equation}
B^{i,j} = \bigoplus_{n_1 + \dots + n_j = i, n_k > 0}
 C^{n_1} \otimes \dots \otimes C^{n_j}.
\end{equation}

So, one grading, let us call it the horizontal one, is just
counting the number of tensor factors, and the other one,
let us call it vertical,
comes from the grading on $C^{>0}(L)$. 
The total grading is their {\em difference}:
$(\Bar C^{>0}(L))^n = \oplus_{i-j=n} B^{i,j}$.
The differential in (\ref{bard}) is decomposed into the 
sum of the vertical differential and the horizontal differential.
The vertical differential depends on the cobracket in $L$,
the horizontal one only depends on the multiplication
in the exterior algebra.

\hfill

Now note that the operation of transition fron $\Bar C^{>0}(L)$
to $\gr_N \Bar C^{>0}(L)$ does not change the horizontal
differential and only changes vertical. The operation
of then taking $\gr_G$ doesn't change the horizontal differential either,
and erases the vertical one.
Moreover, the
complex $\gr_G \gr_N \Bar C^{>0}(L)$ could be identified
with $\gr_G \Bar C^{>0}(L)$, which is just the bar construction
for an exterior coalgebra. Its cohomology
are $\Tor_{C(L)}(\k, \k)$, which could be then
calculated by the Koszul resolution
\cite{PRIDDY}. This calculation shows that its cohomology lie in total degree zero
and are isomorphic to the symmetric coalgebra of $L$.

\hfill

The zeroth cohomology $H^0(\Bar C^{>0}(L))$ could be calculated 
directly. The space of zero-chains of $\Bar C^{>0}(L)$
is the tensor coalgebra generated by $C_1(L) = L$,
that is, the direct sum $\oplus_{i \ge 0} L^{\otimes i}$.
The space
of one-chains is 
$\oplus_{i, j \ge 0} L^{\otimes i} \otimes \Lambda^2L \otimes L^{\otimes j}$.
The differential, from the universal property of tensor coalgebra,
is completely defined by its corestriction to
$\Lambda^2L \subset (\Bar C^{>0}(L))^1$. This is the difference of two maps
$d_1: L \arrow \Lambda^2L$ and $d_2: L \otimes L \arrow \Lambda^2L$,
the first of which is the cobracket in $L$ and the second is 
the alternation map.

\hfill

To make everything more familiar, it is better to take the dual spaces,
especially since we defined the conilpotent coenveloping coalgebra
as the coalgebra topologically dual to the completion of the usual
enveloping algebra in the augmentation ideal. If $\g := L^\vee$ then
 the coalgebra
$H^0(\Bar C^{>0}(L))$ is dual to the algebra $\prod_i \g^{\otimes i}$
factored by the topological ideal generated by elements
of the form $x \otimes y - y \otimes x - [x,y]$, that is, the usual 
(completed) universal enveloping algebra. The theorem is proved.

\hfill

\remark
Positselski's proof of usual PBW theorem uses the
Chevalley-Eilenberg {\em coalgebra} $C_*(\g)$
of a Lie algebra $\g$. We include it here 
for the convenience of the reader, since
the proof in \cite{POSIC} is about more general situation.
It goes as follows. First,
consider (the coaugmented variant of)
the cobar construction $\Cobar C_*(\g)$.
It is a dg-algebra, whose zero cohomology, as we've just seen,
are isomorphic to the usual (uncompleted) universal enveloping algebra 
$U\g$. Consider the filtration on $\Cobar C_*(\g)$ induced from 
the "stupid" filtration on $C_*(\g)$ given by 
$G_iC_*(\g) := \oplus_{k=-i}^0 C_k(\g)$. It is positive ascending 
exhaustive filtration on $\Cobar C_*(\g)$, whose
associated graded could be identified with $\Cobar C_*(V)$,
where $V$ is an abelian Lie algebra which is isomorphic to $\g$
as a vector space. In this case $C_*(V)$ is just
an exterior coalgebra without any differential. As the cobar construction
for a coalgebra $C$ without differential calculates 
$\Ext_C(\k, \k)$, the Ext algebra from the trivial $C$-comodule
to itself, the cohomology algebra $H(\Cobar C_*(V))$
 is the symmetric algebra
$\Sym(V)$ in degree zero and zero in higher degrees, as could be seen
from the same Koszul resolution. Since the filtration $G$
is exhaustive, the associated spectral sequence converges and,
moreover, degenerates at the first page by the dimension reasons.
It follows that 

\begin{equation}
\gr_FH^0(\Cobar C_*(\g)) = \gr_F U\g = 
H^0(\gr_F \Cobar C_*(\g)) = \Sym \g.
\end{equation}

This is the statement of the classical PBW theorem. The difference
between the Lie algebras and Lie coalgebras is that the filtration
on the cobar construction induced from the "stupid" filtration
on the Chevalley-Eilenberg coalgebra
is ascending and exhaustive (in other
words, cocomplete), but the analogous filtration on the bar construction
of the Chevalley-Eilenberg algebra
is descending and not complete. To safely transit from the associated graded
complex to the original one,
in the Lie coalgebra case one needs an extra filtration, which 
is provided by the conilpotency condition.

\hfill

\theorem
Let $L$ be a conilpotent Lie coalgebra. Let $H(L)$ be a minimal $\ainf$-algebra
quasiisomorphic to $C(L)$. Then $H$ is 1-generated.

\hfill

{\bf Proof:} The vanishing of $F^1H^{\ge 1}(\Bar H(L))$ is tautological,
since it is a subspace of $H^{\ge 1}(\Bar H(L))$ which, by \ref{ssequence},
is isomorphic to $H^{\ge 1}(\Bar C(L))$, which vanishes by \ref{barcobar}.

\hfill

\section{C-infinity algebras}\label{OPERADS}

If $A$ is a commutative dg-algebra,
then $\Bar A$ is not only a dg-coalgebra, but also
a dg-bialgebra: the bar-differential
derives the shuffle multiplication 
$\black: \Bar A \otimes \Bar A \arrow \Bar A$. 
If $I$ is the augmentation ideal of $\Bar A$, then
the {\bf complex of indecomposables} $I/I \black I$
obtains the structure of a dg-Lie coalgebra.
This dg-Lie coalgebra is called the {\bf Harrison complex}
of $A$ and will be denoted by
$\Har A$. By the dual Milnor-Moore theorem 
(see \cite{BARR}),
the dg-bialgebra $\Bar A$ is the conilpotent 
 coenveloping coalgebra of $\Har A$.

\hfill

\remark
An important consequence of the dual Milnor-Moore theorem
is that in characteristic zero
$\Har A$ is the {\em direct summand} of $\Bar A$.
In particular,
the embedding $\Har A \arrow \Bar A$ induces an injection on cohomology.

\hfill

\definition
An $\ainf$-algebra $A$ is called a {\bf C-infinity}
algebra, if the coderivation $\delta$ of $\Bar A$
derives the shuffle product $\black$. Alternatively,
a $\cinf$-algebra is a graded vector space $A$
endowed with a square-zero coderivation of the
free Lie coalgebra generated by $A[1]$.
A $\cinf$-morphism between two $\ainf$-algebras
is an $\ainf$-morphism that additionally preserve
shuffle multiplication. Alternatively,
this is a morphism of the corresponding
 Harrison dg-Lie coalgebras.

\hfill

A $\cinf$-algebra is minimal if it is minimal as an $\ainf$-algebra.
The minimal model of an $\ainf$-algebra
admits a $\cinf$-structure, and an $\ainf$-quasiisomorphism
between a $\cinf$-algebra and its minimal model could be chosen
to be a $\cinf$-quasiisomorphism. See, e.g. \cite{GETZLER}.

\hfill

The proof of the \ref{barcobar}, {\em mutatis mutandis}, shows that
the cohomology of $\Har C(L)$ for the conilpotent Lie coalgebra
$C$ vanish in positive degrees and are isomorphic to $C$ in degree zero.
This, of course, is a more general fact which we, for the sake of
brevity, will only state there;
the proofs are largely analogous to the proofs in the Section \ref{PBW}.
For the definitions, one is referred to \cite{LODAYVALETTE}.

\hfill

\theorem
Let $\mathcal{P}$ and $\mathcal{Q}$ be a pair of Koszul dual
operads. Let $C$ be a conilpotent coalgebra over $\mathcal{P}$.
Then the minimal model of its Koszul dual $\mathcal{Q}$-algebra
$\Cobar_\mathcal{P} C$ is 1-generated.

\noindent {\sc {\bf G.P.}:
Northwestern University, \\
Evanston, IL, USA \\
also: \\
Higher School of Economics, \\
Moscow, Russia \\
\tt  grigorypapayanov2020@u.northwestern.edu}

\end{document}